\documentstyle[12pt]{article}

\newcommand{\ellip}{{\cal E}_m}

\newtheorem{lemma}{{\sc Lemma}}
\newtheorem{thm}{{\sc Theorem}}
\newtheorem{proposition}{{\sc Proposition}} 

\newcommand{\Comp}{{\mathbf{C}}}
\newcommand{\Natl}{{\mathbf{N}}}
\newcommand{\Real}{{\mathbf{R}}}

\makeatletter 
 
\@addtoreset{equation}{section}
\makeatother

\newcommand{\bysame}{%
  \leavevmode\hbox to 3em{\hrulefill}\,}

\begin{document}

\title{Singularities of the Bergman kernel
 \\ for certain weakly pseudoconvex domains}
\author{By {\it Joe Kamimoto} \\
{\small Graduate School of Mathematical Sciences,  
The University of Tokyo,} \\
{\small 3-8-1, Komaba, Meguro, 
Tokyo, 153 Japan.}\\
{\it E-mail} : {\tt kamimoto@ms.u-tokyo.ac.jp} }
\date{}
\maketitle



\begin{abstract}
\footnote{{\it Math Subject Classification.} 32A40, 32F15, 32H10.} 

Consider the Bergman kernel $K^B(z)$ of the domain 
$\ellip = \{z \in \Comp^n 
; \sum_{j=1}^n |z_j|^{2m_j}<1 \}$, 
where $m=(m_1,\ldots,m_n) \in \Natl^n$ 
and $m_n \neq 1$.  
Let $z^0 \in \partial \ellip$ 
be any weakly pseudoconvex point, 
$k \in \Natl$ the degenerate rank of the Levi form 
at $z^0$. 
An explicit formula for $K^B(z)$ modulo analytic functions 
is given in terms of the polar coordinates 
$(t_1, \ldots, t_k, r)$ around $z^0$. 
This formula provides detailed information 
about the singularities of $K^B(z)$, 
which improves the result of A. Bonami and 
N. Lohou\'e \cite{bol}. 
A similar result is established also 
for the Szeg\"o kernel $K^S(z)$ of $\ellip$. 
\end{abstract}

\clearpage

\section{Introduction and main result}

Let $\Omega$ be a bounded domain with smooth boundary 
in ${\bf C}^n$, 
$B(\Omega)$ the set of holomorphic $L^2$-functions 
on $\Omega$. 
It is well-known that $B(\Omega)$ is a closed linear 
subspace of the Hilbert space $L^2(\Omega)$. 
The {\it Bergman kernel} $K^B(z)$ of the 
domain $\Omega$ is defined by 
$$
K^B(z) 
=\sum_j |\phi_j(z)|^2, 
$$ 
where $\{\phi_j\}$ is a complete orthonormal basis for 
$B(\Omega)$. 
The above series converges uniformly on any compact 
subset of $\Omega$.  
It is very important to investigate the singularities of 
$K^B(z)$.  
This is mainly because they contain much 
information about the analytic and geometric 
invariants of the domain~$\Omega$. 


First we consider the case where $\Omega$ is a 
strongly pseudoconvex domain. 
In this case C. Fefferman \cite{fef}, 
L. Boutet de Monvel and J. Sj\"ostrand \cite{bos} 
obtained the following asymptotic expansion for $K^B(z)$: 
\begin{equation}
K^B(z)=\frac{\varphi^B(z)}{r(z)^{n+1}}+\psi^B(z)\log r(z),\label{eqn:fef}
\end{equation}
where $r$ is a defining function of~$\Omega$, i.e.,
$\Omega=\{z\in\Comp^n;
r(z)>0\}$ and ${\mathrm{grad}}\,r(z)\ne0$ on
$\partial\Omega$. The functions $\varphi^B(z)$ and $\psi^B(z)$
can be expressed as a power series of~$r$. From the viewpoint of
ordinary differential equations, this result may be interpreted
that the Bergman kernel of a strongly pseudoconvex domain has the
singularities of \emph{regular singular type}.

Next we proceed to the case of weakly pseudoconvex 
domain of finite type  
(in the sense of  J. J. Kohn \cite{koh} 
or J. P. D'Angelo \cite{dan2}).  
In this case there is no such strong general result 
that is comparable with (\ref{eqn:fef}) 
in the strongly pseudoconvex case ; 
yet there are many detailed results for 
the Bergman kernels of specific domains. 
We refer to 
\cite{ber},\cite{dan1},\cite{grs},\cite{dan3},\cite{frh},\cite{has}  
for explicit computations, 
to 
\cite{her1},\cite{ohs1},\cite{dho},\cite{cat},\cite{her2},\cite{her3}  
for  estimates of the size and to 
\cite{bsy} for  boundary limits on nontangential cone.  
Especially D. Catlin \cite{cat} and G. Herbort \cite{her3} 
gave  precise estimates of $K^B(z)$ from above and below 
for certain class of domains whose degenerate 
rank of the Levi form equals one. 
In general, 
however, 
the singularities of $K^B(z)$ are so complicated 
that a unified treatment of them seems to 
be difficult. 



In this paper, 
we pick up the specific domains 
\begin{equation}
\ellip = \left\{
z=(z_1,\ldots,z_n) \in \Comp^n \,;\,
 \sum_{j=1}^n |z_j|^{2m_j} < 1 
\right\}, 
\label{eqn:ell}
\end{equation}
where $m=(m_1,\ldots,m_n) \in \Natl^n$ and $m_n \neq 1$, 
to clarify what is happening for the weakly 
pseudoconvex domains of finite type. 
Since $\ellip$ 
is a Reinhardt domain, 
the set of (normalized) monomials 
forms a complete orthonormal basis for $B(\ellip)$. 
Hence $K^B(z)$ can be represented by a convergent power 
series of $(|z_1|^2,\ldots,|z_n|^2)$, 
whose coefficients were explicitly computed in 
\cite{ise},\cite{dan1},\cite{bol}.      
A. Bonami and N. Lohou\'e \cite{bol} 
gave an important  integral representation for 
the Bergman kernel $K^B(z)$ of $\ellip$ (see (2.1)).  
From this representation they 
deduced a detailed information about the singularities 
of $K^B(z)$, though their result is yet to be improved.

From our point of view, 
we briefly review the result of \cite{bol}. 
Let $z^0=(z_1^0, \ldots, z_n^0) \in \partial \ellip$ 
be any boundary point of $\ellip$, 
$k\in {\bf Z}_{\geq 0}$ 
the degenerate rank of the Levi form at $z^0$. 
We say that $z^0$ is a {\it strongly} 
(resp. {\it weakly}) {\it pseudoconvex point} if $k=0$ 
(resp. if $k>0$). 
Let $I,P$ and $Q$ be the subsets of $N=\{1,\ldots,n\}$ 
defined by 
\begin{equation}
\left\{ 
\begin{array}{rl} 
I &= \{ j \in N ; m_j=1 \}, \\ 
P &= \{ j \in N ; z_j^0=0 \} \setminus I,  \\ 
Q &= \{ j \in N ; z_j^0 \neq 0 \} \cup I.  
\end{array}\right.
\label{eqn:ipq}
\end{equation}
Then the degenerate rank $k$ equals the cardinality 
$|P|$ of $P$. 
One of the main results in \cite{bol} (p.181) 
states that the restriction of $K^B(z)$ 
to the subset $V=\{z_j=0;j\in P\}$ 
admits the following expression around $z^0$: 
\begin{equation}
K^{B}(z)
= C_{P}^{B} \ \ \frac{ \prod_{j \in Q} m_{j}^{2} |z_{j}|^{2m_{j}-2}}
  { \{ 1- \sum_{j \in Q}
 |z_{j}|^{2m_{j}} \} ^{|Q|+| \frac{1}{m} |_{P}+1}}
\ \  +  \ \ O(1), 
\label{eqn:o1}
\end{equation}
where $C_P^B$ is a positive constant and 
$|\frac{1}{m}|_P=\sum_{j\in P}\frac{1}{m_j}$. 
The formula (\ref{eqn:o1}) is quite explicit, 
but still weak in the sense that it is valid only 
on the thin set $V$, 
and that the error term $O(1)$ is somewhat too loose. 

Besides \cite{bol} there  are some studies 
on the Bergman kernel (or Szeg\"o kernel) 
of the domain $\ellip$ 
(\cite{ber},\cite{ise},\cite{dan1},\cite{geb},\cite{goz}). 
In the case $m=(1,\ldots,1,m)$, 
explicit expressions for $K^B(z)$ are 
obtained (\cite{ber},\cite{dan1},\cite{bol}, 
see also Remark 2, \S2.3), 
while there seems to be no explicit one for general $m$. 
The recent studies of Gebelt \cite{geb} and 
Gong and Zheng \cite{goz} are very interesting. 
N. W. Gebelt \cite{geb} generalized  
the method of producing 
the asymptotic expansion 
(\ref{eqn:fef}) due to Fefferman  \cite{fef} 
to the weakly pseudoconvex case of $\ellip$ and 
obtained the analogous results about $K^B(z)$ of 
$\ellip$ $(m=(1,\ldots,1,m))$. 
S. Gong and X. Zheng \cite{goz} 
gave a global estimate of $K^B(z)$ 
from above and below.


Now we state our main results. 
Our essential idea is to introduce the new variables 
$(t,r)$, 
which we call the {\it polar coordinates} 
around $z^0$. 
Here  $t=(t_j)_{j\in P}$ is defined by 
$$
t_j(z)^{2m_j}=
\frac{|z_j|^{2m_j}}
{1-\sum_{j \in Q} |z_j|^{2m_j}}      
\,\,\,\,\,\,\,\,\,\,\,
(j \in P),
$$
and $r$ is the defining function of $\ellip$, i.e., 
$$
r(z)
=1-\sum_{j=1}^n |z_j|^{2m_j}. 
$$
We call 
$t$ the {\it angular variables} and 
$r$ the {\it radial variable}, respectively. 
Then the map $F: z \mapsto (t,r)$ 
takes $\ellip$ onto the region: 
$$
D=
\left\{ 
(t,r) \in  {\bf R}^{|P|} \times (0,1] \,;\, 
t_j \geq 0, \sum_{j\in P} t_j^{2m_j} \leq 1-r 
\right\}. 
$$
The accumulation points of $F(z)$ as 
$\ellip \ni z \to z^0$ are precisely those 
points which belong to the set 
$\{0\} \times \overline{\Delta}$, 
where 
$\overline{\Delta}$ is the closure of the 
{\it locally closed} simplex: 
$$
\Delta=
\left\{ 
t=(t_j)_{j \in P} \,;\, 
t_j \geq 0, \sum_{j\in P} t_j^{2m_j} < 1 
\right\}. 
$$

Let $G=U\cap K$ be a locally closed subset of 
an Euclidean space, 
where $U$ is open and $K$ is closed, 
respectively. 
Then  we say that $f\in C^{\omega}(G)$ if 
$f$ is a real analytic function on some open neighborhood 
$V$ of $G$ in $U$, 
where $V$ may depend on $f$. 


The following theorem asserts that the 
asymptotic behavior of $K^B$ as 
$\ellip \ni  z \to z^0$ can be expressed 
most conveniently in terms of the polar coordinates 
$(t,r)$. 

\begin{thm}

There is a function  
$\Phi^B(t) \in C^{\omega}(\Delta)$ such that 
\begin{equation}
 K^{B}(z) 
 \equiv 
\frac{n!}{{\pi}^{n}} 
 \prod_{j \in Q} m_{j}^{2} |z_{j}|^{2m_{j}-2} 
 \frac{ \Phi^{B} (t(z))}
 { {r(z)}^{|Q|+|\frac{1}{m}|_{P}+1}} \ \ \ \ \
 {\rm modulo} \ \ {\rm C}^{\omega} (\{ z^{0} \}). 
\label{eqn:thm1} 
\end{equation}
Here $\Phi^B(t)$ satisfies {\rm (i)} or {\rm (ii)}. 
 
{\rm (i)} \,\,\, 
If $z^0$ is a strongly pseudoconvex point 
$($i.e. $P=\emptyset)$, 
$\Phi^B(t)=1$ identically. 

{\rm (ii)} \,\,\, 
If $z^0$ is a weakly pseudoconvex point 
$($ i.e. $P\neq \emptyset)$, 
then $\Phi^B(t)$ is positive on $\Delta$ 
and is unbounded as $t \in \Delta$ 
approaches  
$\overline{\Delta} \setminus \Delta$.  

\end{thm}

We remark that 
$\prod_{j \in Q} m_j^2 |z_j|^{2m_j-2}$ 
does not contribute the singularities of 
$K^B(z)$ seriously since it is positive near $z^0$. 
Later we shall see  
that  $\Phi^B(t)$ is essentially the Laplace transform of 
a certain auxiliary function expressible 
in terms of Mittag-Leffler's function 
(see (\ref{eqn:232})). 

We mention a few implications of the formula (\ref{eqn:thm1}) 
in order to compare it with 
the known results stated previously.  
First, 
if $z^0$ is a strongly pseudoconvex point, 
i.e.\ $P=\emptyset$, 
then the angular variables $t$ do not appear 
and $\Phi^B(t)=1$ identically, 
and therefore (\ref{eqn:thm1}) 
reproduces the asymptotic expansion (\ref{eqn:fef}) 
due to C. Fefferman \cite{fef}, L. Boutet de Monvel and 
J. Sj\"ostrand \cite{bos}. 
We remark that the logarithmic term in (\ref{eqn:fef}) 
does not appear in the present case. 
Secondly, 
the restriction of (\ref{eqn:thm1}) to 
the subset $V$ is just the substitution 
$t_j(z)=0$ $(j \in P)$  into (\ref{eqn:thm1}), 
which induces the formula (\ref{eqn:o1}) with the error term 
$O(1)$ replaced by a real analytic function. 
Thus the formula (\ref{eqn:thm1}) 
improves that of Bonami and Lohou\'e 
\cite{bol} in the sense that it is valid 
in a wider domain 
and that the error term is more accurate.  

From our theorem, 
we consider the behavior of $K^B(z)$ at a weakly pseudoconvex 
point from the following three angles: 
(a) {\it estimate}, 
(b) {\it boundary limit} and 
(c) {\it asymptotic formula}.  
We assume $z^0$ is a weakly pseudoconvex point and 
define the region ${\cal U}_{\alpha}(z^0)(={\cal U}_{\alpha}) \subset 
\ellip$ 
by 
$$
{\cal U}_{\alpha}=
\left\{
z \in \ellip ; \sum_{j\in P} t_j(z)=\frac{\sum_{j\in P}|z_j|^{2m_j}}
{1-\sum_{j\in Q}|z_j|^{2m_j}} < \frac{1}{\alpha} 
\right\}\,\,\,(\alpha>1).
$$
(a) By the boundedness of $\Phi^B(t)$ in (ii) 
we can precisely estimate the size of $K^B(z)$ on ${\cal U}_{\alpha}$. 
The region ${\cal U}_{\alpha}$ 
reminds us of the {\it admissible approach regions} 
considered in 
\cite{nsw1},\cite{nsw2},\cite{kra2},\cite{kra3},\cite{ala} etc. 
\,(b) The boundary limit of 
$K^B(z) \cdot r(z)^{|Q|+|\frac{1}{m}|_P+1}$ 
as $z \to z^0$ on each ${\cal U}_{\alpha}$ 
is not determined uniquely but  
depends on the angular variables $t$.   
Note that this boundary limit is uniquely determined  
on any nontangential cone. 
\,(c) In view of  (\ref{eqn:thm1}) the polar coordinates $(t,r)$ is 
necessary to understand the asymptotic formula of $K^B(z)$ 
at $z^0$. 
This fact may be interpreted that 
the Bergman kernel has a singularity of 
irregular singular type at a weakly pseudoconvex point. 
The degeneration from the strong pseudoconvexity 
to the weak pseudoconvexity corresponds to 
the process of confluence from the regular singularity 
to the irregular singularity (\cite{mal},\cite{maj}).

In more detail 
we investigate  the structure of singularities of 
the Bergman kernel of $\ellip$. 
The singularities of $\Phi^B(t)$ at 
$\overline{\Delta} \setminus \Delta$  can also be 
expressed in a  form similar to (\ref{eqn:thm1}) 
by introducing  new polar coordinates on the simplex $\Delta$. 
Through the  finite recursive procedure of this type 
we can completely understand the structure of the singularities 
of $K^B(z)$. 
This situation will be  explained more precisely in Section 3. 


  
This paper is organized as follows. 
In Section 2  
we give the proof of Theorem 1. 
We divide the proof into two parts. 
In the first part  
we refine the error term $O(1)$ in (\ref{eqn:o1}). 
In the second part  
we introduce the polar coordinates and 
express the singularities of $K^B(z)$ 
explicitly. 
In Section 3 
we completely investigate the structure of the singularities 
of $\Phi^B(t)$ 
through the finite recursive procedure described above. 
In Section 4  
we give the proof of Lemma~2, 
which is necessary for the proof of Theorem 1. 
In Section 5  
a similar result about the Szeg\"o kernel of $\ellip$ 
is established.


{\it Acknowledgment.}  \,\,
I would like to express my deepest gratitude 
to Katsunori Iwasaki for very useful 
conversations and his kind help 
during the preparation of this paper. 

\section{Proof of Theorem 1}

In this section, 
we give the proof of Theorem 1.        
We write 
$A \stackrel{>}{\sim} B$ to show that 
$|B/A|$ is bounded 
and  when 
$A \stackrel{>}{\sim} B$ 
and $A \stackrel{<}{\sim} B$, 
we write $A \approx B$.


\subsection{Integral representation of Bonami and Lohou\'e}

Bonami and Lohou\'e \cite{bol} give  
the following integral  
representation of the Bergman kernel of 
$\ellip$: 
\begin{equation}
 K^{B} (z)
 =\frac{1}{\pi^{n}}  
 \int_0^{\infty} e^{-\tau} 
 \prod_{j=1}^n F_{m_{j}}(|z_{j}|^{2} \tau^{\frac{1}{m_{j}}} )
 \tau^{|\frac{1}{m}|_{N}} d \tau,      
\label{eqn:blint}
\end{equation}
with 
$$
 F_{m}(u) 
 = m \sum_{\nu =0}^{\infty} 
 \frac{u^{\nu}}{\Gamma( \frac{\nu}{m} + \frac{1}{m} )},    
$$
where $m\in \Natl$.  
Here $F_m$ is the  derivative of Mittag-Leffler's function: 
$$
 E_{m} (u)
 = \sum_{ \nu =0}^{\infty} 
 \frac{u^{\nu}}{\Gamma(\frac{\nu}{m}+1)},    
$$
(i.e.\ $E_{m}' = F_m$). 

We briefly explain the method of Bonami and Lohou\'e 
for obtaining the integral representation of $K^B(z)$. 
As mentioned in the Introduction, 
the following power series representation of $K^B(z)$ 
is given in \cite{ise},\cite{dan1},\cite{bol}: 
\begin{equation}
K^B(z)
=c\sum_{\nu}
\frac{
\Gamma(\sum_{j=1}^n \frac{\nu_j}{m_j} +\sum_{j=1}^n \frac{1}{m_j}+1)
     }{ 
\prod_{j=1}^n \Gamma(\frac{\nu_j}{m_j}+\frac{1}{m_j})
     }
\prod_{j=1}^n|z_j|^{2\nu_j}, 
\label{eqn:psb}
\end{equation}
where $c=\frac{1}{2\pi^n}\prod_{j=1}^n m_j$.  
Next we represent the Gamma function in the numerator 
in terms of the integral expression and 
change the order of the integral and the sum. 
Finally we put in order the sum in the integral, 
then we can obtain (\ref{eqn:blint}).


\subsection{Refinement of the error term $O(1)$}

In this subsection, 
we investigate the error term $O(1)$ in (\ref{eqn:o1}) 
(\cite{bol}, p.181)   
more precisely. 
The argument below 
is almost similar to that of Bonami and Lohou\'e \cite{bol}. 

Throughout this section, 
we investigate the Bergman kernel in a small neighborhood 
of the fixed boundary point 
$z^0=(z_1^0,\ldots,z_n^0) \in \partial \ellip$.  
Let $N,\, I, \, P$ and $Q$ be defined by (\ref{eqn:ipq}). 


The following properties of Mittag-Leffler's function 
are necessary for the computation below. 

\begin{lemma}[\cite{mit},\cite{val},\cite{bol}] 
Regarding  $F_m(u)$ as an entire function 
on the complex plane, 
$F_m(u)$ is expressed in the following form: 
\begin{equation}
F_{m} (u) 
= m^{2} u^{m-1} \chi_{{\cal H}}(u) \, e^{u^{m}} + f_{m}(u),  
\label{eqn:fm}
\end{equation} 
where \,
$ \chi_{\cal H}(u) = \left\{
                \begin{array}{cl}
                          1 &  for \,\,\,\,   u \in {\cal H}
                              :=\{|\arg u|<\frac{\pi}{2m}\} \\
                          0 &   otherwise,   
                \end{array}
                           \right.$  
and the function $f_m(u)$ has the following properties:  
 {\rm (i)}  $f_{m}(u)$ is bounded in ${\bf C}$, 
     {\rm (ii)}  $f_{m}(u)$  is  holomorphic  on  ${\cal H}$,  
and  {\rm (iii)}  there is a positive constant $c$ 
             such that
            $ f_{m}(u) > c > 0$  for $u > 0$  
             and  
             $\lim_{{u \to 0} \atop {u \in {\cal H}}}
             f_{m}(u) = \frac{m}{ \Gamma( \frac{1}{m} ) }$.   
\end{lemma} 

Substituting (\ref{eqn:fm}) 
into the integral representation (\ref{eqn:blint}), 
we have 
$$
K^{B}(z) 
 = \frac{n!}{\pi^n} 
 \sum_{I \subseteq K \subseteq N} I_{K}(z),  
$$
where
\begin{eqnarray}
  \lefteqn{
 I_{K}(z) 
  = \frac{1}{n!} 
        \prod_{j \in K} m_{j}^2 |z_{j}|^{2m_{j}-2} 
   }    \nonumber\\
 && \times \int_0^{\infty}
        e^{- [1 - \sum_{j \in K} |z_{j}|^{2m_{j}}] \tau} 
        \prod_{j \in N \setminus K} f_{m_{j}}(|z_{j}|^{2} \tau^{m_{j}} )
        \tau^{|K| + |\frac{1}{m}|_{N \setminus K} } d \tau.  
\label{eqn:Ik}
\end{eqnarray}
Applying Lemma~1 to (\ref{eqn:Ik}), 
we obtain the following estimate for  $I_K(z)$:  
\begin{equation}               
 I_{K}(z) \approx
           \frac{ \prod_{j \in K} m_{j}^{2}
           |z_{j}|^{2m_{j} - 2} }
           { \left[ 1 - \sum_{j \in K} 
           |z_{j}|^{2m_{j}} \right]
           ^{|K| + |\frac{1}{m}|_{N \setminus K} + 1} } 
           \,\,\,\, {\rm near} \,\,\,  z^{0}. 
           \label{eqn:Ikest}        
\end{equation}  
By (\ref{eqn:Ikest}), 
we know that  $I_K(z)$ is unbounded near $z^0$ 
,if and only if  
$K \supseteq Q$.  
More precisely we have 
\begin{lemma}
$$   
\sum_{K \not\supseteq Q}
I_K(z) \in {\rm C}^{\omega}(\{ z^{0} \}). 
$$
\end{lemma} 
This lemma will be established in Section 4. 
It implies   
\begin{equation}
 K^{B}(z) \ \ \equiv \ \
\frac{n!}{\pi^n} \sum_{K \supseteq Q} I_{K}(z) 
\ \ \ \ {\rm modulo} \ \ C^{\omega} (\{ z^{0} \}).  
\label{eqn:201}
\end{equation}
Now, restricting $K^B(z)$ 
to the set $V=\{z_j=0;j\in P\}$, 
we have 
\begin{equation}
 K^{B}(z)
\equiv C_{P}^{B} \, 
\frac{ \prod_{j \in Q} m_{j}^{2} |z_{j}|^{2m_{j}-2}}
  { \{ 1- \sum_{j \in Q}
 |z_{j}|^{2m_{j}} \} ^{|Q|+| \frac{1}{m} |_{P}+1}}
\ \  {\rm modulo}  \ \ {C}^{\omega}(\{z^0\}).  
\label{eqn:202}
\end{equation}
In fact, 
$I_K(z)$ ($K\neq Q$) vanishes  identically  
on $V$  and  
$f_m(0)=\frac{m}{\Gamma(\frac{1}{m})}$ 
by Lemma~1.  
The above formula is an improvement of (\ref{eqn:o1})  
and if $z^0 \in \partial \ellip$  
is a strongly pseudoconvex  point (i.e.\ $Q=N$), 
then we obtain (i) in the theorem. 


Now  
we suppose that $z^0\in \partial\ellip$ is 
a weakly pseudoconvex point 
and  investigate the behavior of $K^B(z)$ at $z^0$ 
without the above restriction. 
By (\ref{eqn:Ikest}) and (\ref{eqn:201}), 
we obtain a precise estimate from above and below: 
$$
 K^B(z)   
  \approx \sum_{K \supseteq Q}
  \frac{ \prod_{j \in K} m_{j}^{2}
                 |z_{j}|^{2m_{j} - 2} }
                 { r(z)^{|K| + |\frac{1}{m}|_{N \setminus K} + 1} } 
    \ \ \ \ \   {\rm near} \ \ z^0.  
$$
Note that the fact that 
$f_m(u) \approx 1$  for  $u \geq 0$ \,  
plays an essential role in obtaining  the above estimate. 
Furthermore  we would like to 
investigate the asymptotic behavior 
of $K^B(z)$ at $z^0$. 
For this purpose, 
it is an important problem to obtain 
appropriate information about the function $F_m(u)$.     
Bonami and Lohou\'e (\cite{bol}, pp.177-178) indicate that 
the asymptotic expansion of $K^B(z)$ can be obtained 
by using that of the function $F_m(u)$ at infinity. 
But the result obtained in their manner is difficult to
write clearly and 
the meaning of this expansion seems not to be clear. 
In this paper, 
we assert that Lemma~1 is 
sufficient information about $F_m(u)$ 
to obtain the asymptotic formula of $K^B(z)$. 
Instead of more detailed analysis of $F_m(u)$, 
we introduce another geometric idea, 
which will be mentioned in the next subsection.


\subsection{New coordinates}

In this subsection, 
we continue the argument of the previous subsection 
and complete the proof of the theorem. 

From (\ref{eqn:fm}) and  (\ref{eqn:Ik}), 
we have  
\begin{eqnarray}
     \lefteqn{\sum_{K\supseteq Q}I_K(z)
     =\frac{1}{n!} 
     \prod_{j \in Q} m_{j}^{2} |z_{j}|^{2m_{j}-2}}\nonumber\\
&&  \times  \int_0^{\infty} e^{-\tau [1-\sum_{j \in Q} |z_{j}|^{2m_{j}}]} 
     \prod_{j \in P} F_{m_{j}}(|z_{j}|^{2} \tau^{\frac{1}{m_{j}}} )
     \tau^{|Q|+|\frac{1}{m}|_{P}} d \tau. \nonumber
\end{eqnarray}
Now we introduce new variables $t=(t_j)_{j\in P}$, 
where 
$$
 t_{j}^{2m_j}=\frac{|z_{j}|^{2m_{j}}}
{1-\sum_{j \in Q } |z_{j}|^{2m_{j}}} \,\,\,\,\,\,\,
(j \in P).
$$ 
Then we have 
\begin{equation}
  \sum_{K \supseteq Q} I_{K}(z) = 
 \prod_{j \in Q} m_{j}^{2} |z_{j}|^{2m_{j}-2} 
 \frac{ \Phi^{B} (t(z))}
 { {r(z)}^{|Q|+|\frac{1}{m}|_{P}+1} }, 
\label{eqn:231}
\end{equation}
where 
\begin{equation}
  \Phi^{B} (t) =
 \frac{1}{n!}
 \Bigl[ 1 - \sum_{j \in P} t_{j}^{2m_j} \Bigr]
 ^{|Q|+|\frac{1}{m}|_{P}+1}
 \int_0^{\infty} e^{-s} 
 \prod_{j \in P} F_{m_{j}}
 ( t_{j}^2  s^{\frac{1}{m_{j}}})
 s^{|Q|+|\frac{1}{m}|_{P}} ds.    
\label{eqn:232}
\end{equation}
Substituting (\ref{eqn:231}) into (\ref{eqn:201}), 
we obtain (\ref{eqn:thm1}) in Theorem 1. 
Since 
$\prod_{j \in Q} m_j^2 |z_j|^{2m_j-2}$ 
is positive near $z^0$, 
(\ref{eqn:231}) implies 
that the singularities of $K^B(z)$ is essentially 
expressed in terms of 
the polar coordinates $(t,r)$. 
We show the remaining assertions of Theorem 1. 

First,  
we can obtain that 
$\Phi^B(t)$  is  real analytic 
on the locally closed simplex: 
\begin{equation}
\Delta=
\left\{t = (t_j)_{j \in P} \,;\, 
t_j \geq 0,\, \sum_{j\in P} t_j^{2m_j} <1
\right\}, 
\label{eqn:delta}
\end{equation}
in the same fashion as 
in the proof of Lemma~2.   

Next we can obtain: 
\begin{equation}
\Phi^B (t(z)) =
  \frac{1}{n!} 
  [1-\sum_{j \in P} t_j^{2m_j}]^{|Q|+|\frac{1}{m}|_P +1} 
  \sum_{K \subseteq P} 
  J_K(t(z)) 
\label{eqn:233}
\end{equation} 
where 
\begin{eqnarray*}
\lefteqn{J_K(t)
=\prod_{j \in K} m_j^2 t_j^{2m_j-2}} \nonumber\\ 
& & \times \int_0^{\infty} 
e^{-[1-\sum_{j \in K} t_j^{2m_j}]s} 
\prod_{j \in P \setminus K} 
f_{m_j}(t_j^2 s_j^{\frac{1}{m_j}}) 
s^{|Q|+|K|+|\frac{1}{m}|_{P \setminus K}}
ds, 
\end{eqnarray*}
in the same fashion as in Subsection 2.2.   
Here each $J_K(t)$ has the following estimate: 
\begin{equation} 
J_K(t)
\approx 
\frac{\prod_{j \in K} m_j^2 t_j^{2m_j-2}} 
{[1-\sum_{j \in K} t_j^{2m_j}]
^{|Q|+|K|+|\frac{1}{m}|_{P \setminus K}+1}}. 
\label{eqn:234}
\end{equation}

Now  we claim  
(a) $\Phi^B(t)$ is positive on $\Delta$ 
and 
(b) $\Phi^B(t)$ is unbounded as $t \in \Delta$ 
approaches  $\overline{\Delta} \setminus \Delta$.

(a) : Since $J_{\emptyset}(t) \stackrel{>}{\sim} 1$ 
by (\ref{eqn:234}), 
$\Phi^B(t)$ is positive on $\Delta$ by (\ref{eqn:233}). 
 
(b) : We consider the case where 
$t$ approaches  $t^0=(t_j^0)_{j\in P} 
\in \overline{\Delta} \setminus \Delta$. 
Let $P_{[2]}$, $Q_{[2]}$ be the sets defined by  
$P_{[2]}=\{j \in P; t_j^0=0 \},$  
$Q_{[2]}=\{j \in P; t_j^0 \neq 0 \}$ respectively. 
By (\ref{eqn:234}),  
we have 
$$
[1-\sum_{j \in P} t_j^{2m_j}]^{|Q|+|\frac{1}{m}|_P+1}
J_{Q_{[2]}}(t) 
\approx [1-\sum_{j \in Q_{[2]}} t_j^{2m_j}]
^{-|Q_{[2]}|+|\frac{1}{m}|_{Q_{[2]}}}.
$$ 
Since $Q_{[2]}$ is not empty, 
$[1-\sum_{j \in P} t_j^{2m_j}]^{|Q|+|\frac{1}{m}|_P+1} 
J_{Q_{[2]}}(t)$ 
is unbounded as $t \to t^0$.    
Hence  we obtain (b) by (\ref{eqn:233}).   

This completes the proof of Theorem 1.  

{\it Remarks.} \,\,\,{\it 1.}\,\,   
Since 
$|Q|+|K|+|\frac{1}{m}|_{P \setminus K}+1 
\leq |P|+|Q|+1 
= n+1$ 
and 
$1-\sum_{j \in K}t_j^{2m_j} 
\geq 1-\sum_{j \in P}t_j^{2m_j}$, 
we have 
$$
J_K(t) \stackrel{<}{\sim} 
[1-\sum_{j\in P} t_j^{2m_j}]^{-|P|+|\frac{1}{m}|_P} 
$$ 
on $\Delta$. 
By (\ref{eqn:233}), 
we have 
\begin{equation}
\Phi^B(t) 
\stackrel{<}{\sim} 
[1-\sum_{j \in P} t_j^{2m_j}]^{-|P|+|\frac{1}{m}|_P}
\stackrel{<}{\sim} 
r(z)^{-|P|+|\frac{1}{m}|_P} 
\label{eqn:235}
\end{equation}
The above estimate is optimal. 
In fact, 
we have 
$J_P(t) \approx \prod_{j \in P} t_j^{2m_j-2}
[1-\sum_{j \in P}t_j^{2m_j}]^{-|P|+|\frac{1}{m}|_P}$. 
By (\ref{eqn:231}),(\ref{eqn:235}), 
we have 
$$
K^B(z) \stackrel{<}{\sim} \frac{1}{r(z)^{n+1}}. 
$$


{\it 2.} \,\,  
In the case $m=(1,\ldots,1,m)$, 
we can obtain the following closed expression of $K^B(z)$: 
$$
 K^{B} (z) 
= \frac{n!}{{\pi}^{n}} 
  \frac{ \Phi^{B} (t)}
 { r^{n+\frac{1}{m}} },
\,\,\,\,
{\rm  with} 
\,\,\,\, 
\Phi^{B} (t) 
=mT^{1-\frac{1}{m}}(1-T)^{n+\frac{1}{m}}
 \frac{d^{n}}{dT^{n}}
 \left( \frac{T^{n-1}}{1-T^{\frac{1}{m}}} \right),   
$$
where $T=t^{2m}$.



\section{Recursive formula}

In this section, 
we investigate the structure of the singularities 
of the Bergman kernel of $\ellip$ in more detail. 
From the viewpoint of ordinary differential equations, 
the argument below reminds us of  
the process of step-by-step 
confluence from the regular singularity 
to the irregular singularity  
(See \cite{kht}). 
In this section, 
the results below can be justified in the same fashion 
in Section 2,  
so we omit the detailed proofs of them. 

We remark that $\Phi^B(t)$ defined by (\ref{eqn:232}) 
takes the same form as $K^B(z)$ in (\ref{eqn:blint}). 
So the argument in Section 2 applies to $\Phi^B(t)$ 
in place of $K^B(z)$, 
and the form of the singularities of 
$\Phi^B(t)$  
can be written in the same fashion as 
in Theorem 1. 
Moreover we can completely understand the singularities of 
$\Phi^B(t)$ 
by finitely many recursive process of this kind.   

We precisely explain this process. 
We suppose that $z^0 \in \partial \ellip$ 
is a weakly pseudoconvex point (i.e.\ $P \neq \emptyset$). 
We inductively define the sets 
$P_{[k]}$, $Q_{[k]} \subset N$, 
the variables $t_{[k]}=(t_{[k],j})_{j \in P_{[k]}}$, 
$r_{[k]}$, the simplex $\Delta_{[k]}$ and 
the function $\Phi_{[k]}(t_{[k]})$ on 
$\Delta_{[k]}$ in the following way.

First we set 
$P_{[1]}=P (\neq \emptyset)$, 
$Q_{[1]}=Q$, 
$t_{[1]}=(t_{[1],j})_{j \in P_{[1]}}=(t_j)_{j\in P}$, 
$r_{[1]}=r$, 
$\Delta_{[1]}=\Delta$ and 
$\Phi_{[1]}(t_{[1]})=\Phi^B (t)$. 
Suppose that the sets 
$P_{[k-1]}(\neq \emptyset), 
\,\,Q_{[k-1]} \subset P$ 
are settled, 
then the simplex $\Delta_{[k-1]}$ is defined by 
\[
\Delta_{[k-1]}=\{t_{[k-1]}=(t_{[k-1],j})_{j \in P_{[k-1]}} ; 
t_{[k-1],j} \geq 0, \,\, 
\sum_{j \in P_{[k-1]}} t_{[k-1],j}^{2m_j}<1 \} 
\subset \Real^{|P_{[k-1]}|}.    
\]
When we select a point 
$t_{[k-1]}^0=(t_{[k-1],j}^0)_{j \in P_{[k]}} 
\in \overline{\Delta_{[k-1]}} \setminus \Delta_{[k-1]}$, 
the sets $P_{[k]},\,Q_{[k]} \subset P_{[k-1]}$ 
are determined by  
\[
 \left\{
\begin{array}{rl}
P_{[k]}&= \{j \in P_{[k-1]}; t_{[k-1],j}^0 =0\}, \\
Q_{[k]}&= \{j \in P_{[k-1]}; t_{[k-1],j}^0 \neq 0\}.  
\end{array}\right.
\] 
Furthermore the variables 
$t_{[k]}=(t_{[k],j})_{j \in P_{[k]}},\,\,
r_{[k]}$ are defined by 
\[ 
\left\{ 
\begin{array}{rl}
t_{[k],j}^{2m_j}&=\frac{t_{[k-1],j}^{2m_j}}
{1-\sum_{j \in P_{[k-1]}} t_{[k-1],j}^{2m_j}} \,\,\,\,\,
(j \in P_{[k]}), \\
r_{[k]}&= 1-\sum_{j \in P_{[k-1]}}t_{[k-1],j}^{2m_j}.   
\end{array}\right.
\]
Then we define the function 
$\Phi_{[k]}(t_{[k]})$ 
on the simplex $\Delta_{[k]}$ 
in the following. 
If $P_{[k]}= \emptyset$, 
then 
$\Phi_{[k]}(t_{[k]})=1$ identically. 
If $P_{[k]} \neq \emptyset$, 
then 
\begin{equation}
\Phi_{[k]}(t_{[k]}) 
=\frac{1}{n!} 
[1-\sum_{j \in P_{[k]}} t_{[k],j}^{2m_j}]
^{a_{[k]}+1}
\int_0^{\infty} 
e^{-s} \prod_{j \in P_{[k]}}F_{m_j}(t_{[k],j}^2 s^{\frac{1}{m_j}}) 
s^{a_{[k]}} ds,  
\label{eqn:30}
\end{equation}
where the constant $a_{[k]}$ is defined by 
$
a_{[k]}=
\sum_{j=1}^k |Q_{[j]}|+
|\frac{1}{m}|_P 
-\sum_{j=2}^k |\frac{1}{m}|_{Q_{[j]}}
$.  
In the above inductive process, 
we have 
$$
P=P_{[1]} 
 \mathop{\supset}_{\ne} P_{[2]}
 \mathop{\supset}_{\ne} \cdots 
 \mathop{\supset}_{\ne} P_{[k-1]} 
 \mathop{\supset}_{\ne} P_{[k]}.
$$
So there exists a positive integer $k^0 \leq |P|$ 
such that $P_{[k^0]}= \emptyset$. 
Thus we have defined 
$P_{[k]},Q_{[k]},\, 
t_{[k]},r_{[k]},\,
\Delta_{[k]},\, 
\Phi_{[k]}(t_{[k]})$ for  
$ k= 1, 2, \ldots, k^0$  recursively.  

We remark that $\Phi_{[k]}(t_{[k]})$ in (\ref{eqn:30}) 
takes the same form as $K^B(z)$ in (\ref{eqn:blint}). 
So we obtain the following proposition for 
$\Phi_{[k]}(t_{[k]})$ in the same manner as we have 
obtained Theorem 1 for $K^B(t)$.

\begin{proposition}
Suppose that $1 \leq k \leq k^0-1$. 
The function $\Phi_{[k]}(t_{[k]})$ is a positive 
and real analytic function  
on $\Delta_{[k]}$  
and is unbounded as 
$t_{[k]}\in \Delta_{[k]}$ approaches $t_{[k]}^0 
\in \overline{\Delta_{[k]}} \setminus \Delta_{[k]}$.
Moreover we have the following recursive formula: 
\begin{equation}
\Phi_{[k]}(t_{[k]}) 
\equiv 
\prod_{j \in Q_{[k+1]}} m_j^{2} t_{[k],j}^{2m_j-2}
\frac{\Phi_{[k+1]}(t_{[k+1]})}
{r_{[k+1]}^{|Q_{[k+1]}|-|\frac{1}{m}|_{Q_{[k+1]}}}}   
\,\,\,{\rm modulo}\,\,\, C^{\omega}(\{t_{[k]}^0\}).    
\label{eqn:prop}
\end{equation}  
\end{proposition} 

We remark that the condition 
$k=k^0$ (resp. $1 \leq k \leq k^0-1$)  corresponds to the 
strongly (resp.  the weakly)  pseudoconvex case 
in Theorem 1. 
The formula (\ref{eqn:prop}) recursively reduces 
$\Phi^B(t)$ to $\Phi_{[k^0]}(t_{[k^0]}) \equiv 1$.  
Hence it may be interpreted that 
the above recursive process resolves the degeneration of 
the Levi form in  the study of singularities of 
the Bergman kernel 
in the weakly pseudoconvex case.


\section{Proof of Lemma 2}

By (\ref{eqn:fm}), we obtain   
$$
\sum_{ K \not\supseteq Q} I_{K}(z)
  =\sum_{ I \subseteq J \subset Q }
   \widetilde{I_{J}}(z), 
$$
where  
\begin{eqnarray}
\lefteqn{
       \widetilde{I_{J}}(z) 
       = \frac{1}{\pi^{n}} 
       \prod_{j \in J} m_{j}^2 |z_{j}|^{2m_{j}-2}
       \int_0^{\infty}
       e^{- [1 - \sum_{j \in J} |z_{j}|^{2m_{j}}] \tau}
    } \nonumber\\ 
&  &  \times \prod_{j \in Q \setminus J} 
      f_{m_{j}}(|z_{j}|^{2} \tau^{\frac{1}{m_{j}}} )
      \prod_{j \in P} F_{m_{j}}(|z_{j}|^{2} \tau^{\frac{1}{m_{j}}} )
      \tau^{|J| + |\frac{1}{m}|_{N \setminus J} } d \tau. 
   \nonumber  
\end{eqnarray}
Thus it is sufficient to show that 
\begin{equation}
 \widetilde{I_{J}}(z) \in {\rm C}^{\omega}(\{ z^{0} \}), 
\label{eqn:41} 
\end{equation}
for $I\subseteq J\subset Q$. 

Let $\widehat{I_{J}}(u)$ 
be the function of complex variables 
$u=(u_1, \ldots, u_n) \in {\bf C}^n$ defined by 
\begin{eqnarray}
\lefteqn{
     \widehat{I_{J}}(u)=
     \int_0^{\infty}
      e^{- [1 - \sum_{j \in J} u_{j}^{m_{j}}] \tau} 
} \nonumber\\
& &      \times \prod_{j \in Q \setminus J} 
      f_{m_{j}}(u_{j} \tau^{\frac{1}{m_{j}}} )
      \prod_{j \in P} F_{m_{j}}(u_{j} \tau^{\frac{1}{m_{j}}} )
      \tau^{|J| + |\frac{1}{m}|_{N \setminus J} } d \tau,  
\label{eqn:42}
\end{eqnarray}
In order to obtain (\ref{eqn:41}), 
it is sufficient to show that 
there exists a neighborhood in $\Comp^n$ 
of $u^0=(u_1^0,\ldots,u_n^0)
    :=(|z_{1}^{0}|^{2}, \cdots ,|z_{n}^{0}|^{2})$ 
such that $\widehat{I_J}(u)$ is holomorphic there. 
Note that $\prod_{j\in J} m_j^2 |z_j|^{2m_j-2}$ 
is real analytic at $z^0$. 

Now let ${\cal N}_J$ be the neighborhood of $u^0$ 
defined by 
\begin{eqnarray}
\lefteqn{ {\cal N}_{J} = 
          \Bigl\{u \in {\bf C}^n ;  |u_j| < 1 \,\, 
          {\rm for} \, j \in P, }\nonumber\\
&  &      |u_{j}-u_{j}^{0}| < \frac{u_{j}^{0}}{2} \,
          {\rm for} \, j \in Q \setminus J  \, ,{\rm and} \,
          1- \sum_{j \in J \cup P} |u_{j}|^{m_{j}} > \frac{\delta}{2}          
          \Bigr\}, \nonumber
\end{eqnarray} 
where $\delta = 1 - \sum_{j \in J \cup P} u_{j}^{0} > 0$. 
We show that $\widehat{I_{J}}(u)$ is holomorphic in 
${\cal N}_J$. 
Since $F_m$ is an entire function 
and $f_m$ is holomorphic in 
the sector $\{u;|\arg u|<\frac{\pi}{2m}\}$  
by Lemma~1 (ii), 
the integrand of (\ref{eqn:42}) is holomorphic 
in ${\cal N}_J$ for $\tau >0$. 
Each partial derivative of the integrand in (\ref{eqn:42}) 
with respect to 
$u_j$   
is continuous on ${\cal N}_J \times (0,\infty)$. 
By Lemma~1, 
we have 
$$ 
|F_m(u)|\leq c |u|^{m-1}e^{|u|^m} \,\,{\rm and}\,\, |f_m(u)|\leq c, 
$$
on $\Comp$, 
where $c$ is a positive constant.  
Thus we have  
\begin{eqnarray*}
|\widehat{I_{J}}(u)|  
 &\le& \int_0^{\infty}
       e^{- [1 - \sum_{j \in J} |u_{j}|^{m_{j}}] \tau} \\
 &   & \times \prod_{j \in Q \setminus J}
       |f_{m_{j}}(u_{j} \tau^{\frac{1}{m_{j}}} )| 
       \prod_{j \in P}|F_{m_{j}}(u_{j} \tau^{\frac{1}{m_{j}}} )|
       \: \tau^{|J|+|\frac{1}{m}|_{N \setminus J}} 
       d \tau \\   
 &\le& c^{n+1-|J|} 
       \int_{0}^{\infty}
       e^{-[1-\sum_{j \in J \cup P}|u_j|^{m_j}] \tau}
       \tau^{(m-1)|P|+|J|+|\frac{1}{m}|_{N \setminus J}} 
       d \tau \\
 &\le& c^{n+1 - |J|} 
       \frac{\Gamma((m-1)|P|+|J|+|\frac{1}{m}|_{N \setminus J}+1)}
       {(\delta/2 )
       ^{(m-1)|P|+|J|+|\frac{1}{m}|_{N \setminus J}+1}} 
\end{eqnarray*}
on ${\cal N}_J$. 
Hence we can see that 
$\widehat{I_J}(u)$ is holomorphic on ${\cal N}_J$ 
by the above inequalities.  



{\it Remark}. \,\,\, 
Consider the smoothness of the Bergman kernel 
$K^B(z,w)$ off the diagonal 
(i.e.\ $\Delta:=\{(z,w); z=w \in \partial \ellip\}$). 
Here $K^B(z,w)=\sum_{j} \phi_j(z) \overline{\phi_j(w)}$, 
where $\{\phi_j\}_j$ is as in the Introduction. 
We have 
$$
K^B(z,w) \in \overline{\ellip} \times \overline{\ellip} 
\setminus \Delta. 
$$
This can be obtained by 
putting together the proof of Lemma~2 and 
the argument in \cite{bol}, pp.170-171.


\section{The Szeg\"o kernel of $\ellip$}

In this section, 
we establish  a  result similar to Theorem 1  
for the Szeg\"o kernel of $\ellip$. 
The result below is obtained 
in the same fashion as in the case of 
the Bergman kernel 
and we omit the proof. 


Let $\Omega$ be a bounded domain in $\Comp^n$ 
with smooth boundary. 
Specify a surface element $\sigma$ 
on the boundary $\partial \Omega$, 
and denote by $H_{\sigma}^{2}(\Omega)$ 
the set of holomorphic functions in $\Omega$ 
having $L^2$-boundary values with respect to $\sigma$. 
The {\it Szeg\"o kernel} of $\Omega$ 
(with respect to $\sigma$) 
is defined by 
$$
K^S(z,w)=\sum_{j}
|\tilde{\phi_j}(z)|^2 
$$
where $\{\tilde{\phi_j}\}$ is a complete 
orthonormal basis for $H_{\sigma}^2(\Omega)$. 

We study the Szeg\"o kernel of $\ellip$ 
with respect to 
the surface element 
which is introduced by 
Bonami and Lohou\'e in \cite{bol} 
(they denote the surface element by $d \mu_{\alpha}$).

An integral representation of $K^S(z)$ is 
obtained  in the same  fashion as in the case of 
the Bergman kernel: 
\begin{equation}
  K^{S} (z)
 =\frac{1}{2 \pi^{n}}  
 \int_0^{\infty} e^{-\tau} 
 \prod_{j=1}^n F_{m_{j}}(|z_{j}|^{2} \tau^{\frac{1}{m_{j}}} )
 \tau^{|\frac{1}{m}|_{N}-1} d \tau. 
\label{eqn:51}
\end{equation}

Since the difference 
between (\ref{eqn:blint}) and (\ref{eqn:51}) does not 
give any essential influence on the argument in 
Section 2, we have a similar result about the singularities 
of $K^S(z)$. 


\begin{thm}

There is a function  
$\Phi^S(t) \in C^{\omega}(\Delta)$ such that 
$$
 K^{S}(z) 
 \equiv 
\frac{(n-1)!}{2 {\pi}^{n}} 
 \prod_{j \in Q} m_{j}^{2} |z_{j}|^{2m_{j}-2} 
 \frac{ \Phi^{S} (t(z))}
 { {r(z)}^{|Q|+|\frac{1}{m}|_{P}}} \ \ \ \ \
 {\rm modulo} \ \ {\rm C}^{\omega} (\{ z^{0} \}). 
$$
Here $\Phi^S(t)$ 
also has the same properties as in Theorem 1 for $\Phi^B(t)$ . 
 
\end{thm}
 

{\it Remarks.} \,\,{\it 1.} \,\, 
The precise expression of $\Phi^S(t)$   
is the following: 
$$
\Phi^{S}(t)
 = \frac{1}{(n-1)!}  
   \Bigl[ 1 - \sum_{j \in P} t_{j}^{2m_j} \Bigr]^
   {|Q|+|\frac{1}{m}|_{P}} 
   \int_0^{\infty} e^{-s} 
   \prod_{j \in P}
   F_{m_{j}}( t_{j}^2 s^{\frac{1}{m_{j}}} )
   s^{|Q|+|\frac{1}{m}|_{P}-1} ds.  
$$
In the case $m=(1,\ldots,1,m)$, 
we have the following closed expression:   
$$ 
 K^{S} (z) 
= \frac{(n\!\!-\!\!1)!}{2{\pi}^{n}} 
  \frac{ \Phi^{S} (t)}
 { r^{n-1+\frac{1}{m}} },
$$ 
with 
$$
\Phi^{S} (t) 
=mT^{1-\frac{1}{m}}(1-T)^{n-1+\frac{1}{m}}
 \frac{d^{n-1}}{dT^{n-1}}
 \left( \frac{T^{n-2}}{1-T^{\frac{1}{m}}} \right),  
$$
where $T=t^{2m}$. 

{\it 2.}\,\, 
Consider 
the smoothness of  the Szeg\"o kernel 
$K^S(z,w)(:=\sum_j \tilde{\phi_j}(z) 
\overline{\tilde{\phi_j}(w)})$ off the diagonal. 
We  obtain 
$$
K^S(z,w) \in \overline{\ellip} \times \overline{\ellip} 
\setminus \Delta, 
$$
in the same fashion as in the case of the Bergman kernel. 
See Remark in Section 4.



\end{document}